\documentclass{article}

\usepackage[utf8]{inputenc}
\usepackage{amsmath,amsfonts,amsthm,physics,mathrsfs,bbm,graphicx,enumerate, tikz-cd, amssymb}
\usepackage[margin=1in]{geometry}
\usepackage{xcolor}
\usepackage{cite}

\DeclareMathOperator{\diag}{diag}
\DeclareMathOperator{\Out}{Out}
\DeclareMathOperator{\Aut}{Aut}

\DeclareMathOperator{\R}{\mathbb{R}}
\DeclareMathOperator{\C}{\mathbb{C}}
\DeclareMathOperator{\Z}{\mathbb{Z}}
\DeclareMathOperator{\Hom}{Hom}

\DeclareMathOperator{\End}{End}

\DeclareMathOperator{\PSL}{PSL}
\DeclareMathOperator{\PGL}{PGL}
\DeclareMathOperator{\GL}{GL}
\DeclareMathOperator{\SL}{SL}

\DeclareMathOperator{\Mod}{Mod}
\DeclareMathOperator{\im}{im}

\newcommand{\K}{\mathbb K}

\renewcommand{\frak}{\mathfrak}

\newcommand{\inverse}{^{-1}}

\newcommand{\Ad}{\mathsf {Ad}}

\theoremstyle{plain}
\newtheorem{theorem}{Theorem}[section]
\newtheorem{lemma}[theorem]{Lemma}
\newtheorem{claim}[theorem]{Claim}
\newtheorem{proposition}[theorem]{Proposition}

\newtheorem{remark}[theorem]{Remark}

\theoremstyle{definition}

\theoremstyle{remark}

\graphicspath{{figures/}}

\begin{document}

\title{The Torelli Group and Euler Class Zero $\PSL(2,\R)$ Representations}
\author{Marit Bobb \and James Farre \and Peter Smillie}
\date{}
\maketitle

\begin{abstract}
    We prove that every representation of a closed surface group of genus at least $3$ into $\PSL(2,\R)$ that is sufficiently close (but not equal to) a representation that has non-compact solvable image has a simple closed separating curve that is mapped to a non-trivial elliptic element in $\PSL(2,\R)$. This can be seen as evidence for an affirmative answer to a question of Bowditch and a conjecture of Goldman.  
\end{abstract}
\section{Introduction}

Geometry, topology, and dynamics interact fruitfully in the
study of surfaces, where hyperbolic geometry plays a central
role. Let $S$ be a closed oriented surface of genus $g\ge 2$,
with fundamental group $\pi$. Hyperbolic structures on $S$
are encoded by conjugacy classes of discrete faithful
representations $\pi\to\PSL(2,\R)$. The mapping class group
$\Mod(S)$ acts properly discontinuously on this locus,
with quotient (two copies of) the Riemann moduli space.
We are interested in the dynamics of $\Mod(S)$ on the
larger space of representations, particularly on the
component containing the trivial representation.

We define the \emph{space of characters} to be
\[
\mathcal X \overset{\mathrm{def}}{=}
\Hom(\pi,\PSL(2,\R))/\sim,
\]
with its quotient topology, where $\sim$ is the equivalence
relation generated by $\rho\sim\rho'$ whenever
$\rho'\in\overline{\PSL(2,\R).\rho}$.
The \emph{diagonal singular locus} is
\[
\mathcal Z \overset{\mathrm{def}}{=}
\{[\rho]\in\mathcal X: [\rho] \text{ has a diagonal representative }\}
\]
It is known that $\mathcal X$ is a smooth manifold equipped with a
$\Mod(S)$-invariant symplectic form away from the full singular locus of non-elementary characters \cite{Goldman:symplectic}.

The \emph{Euler number} of a representation
$\rho:\pi\to\PSL(2,\R)$ satisfies the classical
\emph{Milnor--Wood inequality}
\[
|eu(\rho)|\le 2g-2.
\]
It is locally constant and descends to $\mathcal X$.
Goldman proved that the connected components of
$\mathcal X$ are precisely the fibers
$\mathcal X_i=eu\inverse(i)$,
for $-2g+2\le i\le 2g-2$ \cite{Goldman:components}.
The component $\mathcal X_{2g-2}$ identifies with
Fricke space \cite{Goldman:thesis}, whereas
$\mathcal X_0$ contains $\mathcal Z$.
Each component is preserved by $\Mod(S)$.
Goldman conjectured that the action on $\mathcal X_i$
is ergodic whenever $|i|<2g-2$
\cite[Conjecture 3.1]{Goldman:conjecture}.
March\`e and Wolff resolved this conjecture in genus two,
where a special low-genus phenomenon occurs
\cite{MW:genustwo, MW:sex}.

Motivated by his work on representations of the
punctured torus group, Bowditch asked whether every
$\rho\in\Hom(\pi,\PSL(2,\R))$ with
$|eu(\rho)|<2g-2$ admits an essential simple closed
curve $\gamma\subset S$ such that
$|\tr(\rho(\gamma))|\le 2$
\cite[Question C]{Bowditch}.
An affirmative answer for almost every representation
in $\mathcal X_i$ implies Goldman's conjecture on
that component, except when $(g,i)=(2,0)$
\cite[Theorem 1.6]{MW:genustwo}.

Our main theorem provides evidence for Goldman's
conjecture on $\mathcal X_0$ when $g\ge 3$.

\begin{theorem}\label{thm: exists neighborhood R}
For a closed oriented surface $S$ of genus at least $3$,
there is an $\Mod(S)$-invariant open neighborhood $U$ of
$\mathcal Z\setminus\{[id]\}$ in $\mathcal X_0$ 
such that every $[\rho]\in U$ admits a simple closed
separating curve $\eta$ with
$|\tr(\rho(\eta))|\le 2$.
Moreover, if $[\rho]\in U\setminus \mathcal Z$, then
$\eta$ can be chosen so that $|\tr(\rho(\eta))|<2$, i.e., 
$\rho(\eta)$ is elliptic and non-trivial.
\end{theorem}

Morally, Theorem \ref{thm: exists neighborhood R} is proved by studying the infinitesimal action of the \emph{Torelli subgroup} $\mathcal I \le \Mod(S)$ on $\mathcal X$ at $\mathcal Z$. 
Indeed, $\mathcal Z$ can be identified with the locus of diagonal representations $\pi \to \PSL(2,\R)$, up to an order two symmetry exchanging the two entries.  
Thus $\mathcal Z \cong H^1(S,\R)/\pm1$, and the action of $\mathcal I$ fixes $\mathcal Z$ pointwise, inducing a natural tangent action.  In particular, $\mathcal Z$ is a non-compact subvariety.  It is possible to show that $U$ has infinite volume.  

We use work of Ghazouani \cite{Selim:affC,Selim:orbitclosure}  together with the minimality of the $\Mod(S)$ action on lines in $H^1(S,\R)$ and show that there are many directions where the trace function for a separating curve decreases to second order, and use this to conclude our main theorem.

\subsection*{AI disclosure}
The main idea in this paper has been clear to the authors for a few years now, but we struggled to implement a technical part of our proof strategy using algebro-geometric methods. 
We gave the notes we made for this project to ChatGPT, which gave an elementary version of the technical part of our argument that we were stuck on.  We present this as Proposition \ref{prop: zooming in}. 
We are very pleased with this argument, which follows  closely our original insight, and enjoyed the process of making it presentable to a wider audience.  

\section{Deformation theory}

\subsection{Cohomology}

In this section, we recall the deformation theory of representations of surface groups and its relationship to cohomology.  Good references for this material are \cite[\S1.3]{Goldman:conjecture} and \cite{JM:deformations}.

\subsubsection{Tangent cocycles}
Let $G$ be a Lie group and suppose $\rho_t: \pi \to G$ is a smooth path of representations with $\rho_0 = \rho$.
We identify the Lie algebra $\mathfrak g = T_{id}G$ of $G$ and the space of right invariant vector fields on $G$ as usual.
The derivative at time zero, for each $\gamma\in \pi$ is
\[u(\gamma) = \left.\frac{d}{dt}\right|_{t = 0}\rho_t(\gamma)\rho(\gamma)\inverse\in \mathfrak g. \]
It is straightforward to check that
\begin{equation}\label{eqn: cocycle condition}
    u(\gamma_1\gamma_2)  = u(\gamma_1) + \mathsf{Ad}_{\rho(\gamma_1)}u(\gamma_2).
\end{equation}
Let $Z^1(\pi, \mathfrak g_{\mathsf {Ad}_\rho})$ be the space of functions $u: \pi \to \mathfrak g$ satisfying the cocycle condition \eqref{eqn: cocycle condition}. 
The coboundaries $B^1(\pi,\mathfrak g_{\Ad_{\rho}})$ are those of the form 
\[\delta X(\gamma) = (I-\Ad_{\rho(\gamma)})X,\]
for some $X \in \mathfrak g$, and come from trivial deformations (tangent to a path of conjugations $g_t\rho g_t\inverse$).  
The cohomology $H^1(\pi, \mathfrak g_{\Ad_\rho})$ is the quotient of the cocycles by the coboundaries.  

\subsubsection{Cocycles and coboundaries}

More generally, let $V$ be a vector space, and let $\rho: \pi \to \GL(V)$ be a representation.  Then a $k$-cochain is a map
\[c: \pi^k \to V,\]
whereas a $0$-cochain is just an element of $V$.
The coboundary operator is defined by
\[ \delta c (\gamma_0, \ldots,\gamma_k) = -\rho(\gamma_0).c(\gamma_1, \ldots, \gamma_k) + \sum_{i = 0}^{k-1} (-1)^{i}c(\gamma_0,\ldots, \gamma_i\gamma_{i+1}, \ldots,\gamma_k) +(-1)^{k} c(\gamma_0, \ldots, \gamma_{k-1}).\]
It takes $k$-cochains to $(k+1)$-cochains and satisfies $\delta\circ \delta =0$.
Cohomology groups $H^k(\pi, V_\rho)$ are defined as usual as the quotient of cocycles ($\ker\delta$) by  coboundaries ($\im\delta$).  These definitions recover $H^1(\pi,\mathfrak g_{\Ad_\rho})$ from above.

\subsubsection{Operations in cohomology}

Let $\rho : \pi \to \GL(U)$, $\sigma: \pi \to \GL(V)$, and $\tau: \pi \to \GL(W)$ be homomorphisms, and suppose we have a $\pi$-equivariant bilinear pairing
\[\mu: U\times  V \to W,\]
where the equivariance condition is $\mu(\rho(\gamma)u, \sigma(\gamma)v) = \tau(\gamma)\mu(u, v)$ for all $\gamma \in \pi$.

Define the cup product 
\[\cdot \smile \cdot : Z^1(\pi,U_{\rho}) \times Z^1(\pi, V_{\sigma}) \to Z^2(\pi, W_{\tau})\]
by the rule 
\[(u\smile v)(\alpha, \beta) = \mu(u(\alpha), \sigma(\alpha)v(\beta)).\]
One checks that 
\[\delta(u\smile v) = (\delta u)\smile v -u\smile (\delta v),\]
so that the cup product descends to a map
\[\cdot \smile \cdot : H^1(\pi,U_{\rho}) \times H^1(\pi, V_{\sigma}) \to H^2(\pi, W_{\tau})\]
at the level of cohomology.
\medskip

More concretely, we will be interested in the following two coefficient pairings.
\begin{enumerate}
    \item If $\pi$ acts by linear transformations on a vector space $V$, then it also acts by linear transformations on the dual space $V^*$ preserving the canonical evaluation pairing
    \[\langle \cdot, \cdot \rangle : V\times V^* \to \R.\]
    Explicitly, the action on $V^*$ is given by $\gamma.\alpha = \alpha\circ \rho(\gamma)\inverse$, where $\alpha \in \Hom(V,\R)$.
    \item For a Lie group $G$ with Lie algebra $\mathfrak g$ and a homomorphism $\rho: \pi \to G$, the Lie bracket \[[\cdot, \cdot]: \mathfrak g \times \mathfrak g \to \mathfrak g\]
    is $\Ad_\rho$-equivariant.
\end{enumerate}

The first example yields a map 
\[\cdot \smile \cdot : H^1(\pi, V_\rho) \times H^1(\pi, V^*_\rho)\to H^2(\pi, \R)\cong \R,\]
which is anti-symmetric in the sense that 
\[[u]\smile [\alpha] = - [\alpha]\smile [u],\]
for all $u \in Z^1(\pi, V_\rho)$ and $\alpha \in Z^1(\pi, V^*_\rho)$.  

The second example gives a map 
\[[\cdot, \cdot] : H^1(\pi, \mathfrak g_{\Ad_\rho}) \times H^1(\pi, \mathfrak g_{\Ad_\rho}) \to H^2(\pi, \mathfrak g_{\Ad_\rho}),\]
which is  symmetric.

\subsection{Representations to second order}
Although we work mainly with representations into $\SL(2,\R)$, all results in this section remain valid over $\C$.  Fix a standard presentation $F_{2g} \to \pi$, where $F_{2g}$ is the free group on $2g$ generators.

Let $\K = \R$ or $\C$, and suppose $\rho: F_{2g} \to \SL(2,\K[t])$ is a representation.  Say that $\rho$ \emph{factors through $\pi$ to second order} if 
\[\rho\left(\prod_{i =1}^g[\gamma_{2i-1},\gamma_{2i}]\right) = I \mod t^3,\]
or in other words if $\rho$ factors through $\pi$ as a representation into $\SL(2, \K[t]/(t^3))$.
The reduction of $\rho$ $\mod t$ is a representation into $\SL(2,\K)$, which we denote $\rho_0$.

The following material is fairly standard (cf. \cite[\S1.3]{Goldman:conjecture} or \cite[\S4]{JM:deformations}).  We have chosen to give a self-contained proof because we did not find it in the form we need for the proof of Proposition \ref{prop: zooming in}.
\begin{proposition}\label{prop: deformation theory}
     Suppose $\rho: F_{2g} \to \SL(n,\K[t])$ factors through $\pi$ to second order, and write 
    \[\rho(\gamma) = (I +tZ(\gamma) +t^2W(\gamma))\rho_0(\gamma) \mod t^3,\]
    where $Z(\gamma)$ and $W(\gamma) \in M_n(\K)$.
    Then $Z\in Z^1(\pi, \mathfrak g_{\Ad_{\rho_0}})$, and $[Z, Z] =0\in H^2(\pi,\mathfrak g_{\Ad_{\rho_0}})$, where $\mathfrak g$ is the Lie algebra of $\SL(n,\K)$.
\end{proposition}

\begin{proof}
    Recall Jacobi's formula that for a smooth path of matrices $A(t) \in \SL(n, \K)$,  \[(\log \det A(t))' = \tr (A'(t)A(t)\inverse).\]   Since $\rho(\gamma) \in \SL(n,\K[t])$, $\log\det\rho(\gamma) = 0$, so expanding  about $t =0$ and comparing coefficients of $t$ gives
    \begin{equation}\label{eqn: traceless 1}
         \tr Z = 0
    \end{equation}    
    while comparing coefficients of $t^2$ gives
    \begin{equation}\label{eqn: traceless 2}
        \tr (2W - Z^2) = 0.
    \end{equation}
    
    Since $\rho$ factors through $\pi$ to second order,
    $\rho(\alpha\beta) = \rho(\alpha)\rho(\beta)\mod t^3$.  Multiplying in $\SL(n,\K[t]/(t^3))$ and comparing the coefficients gives 
    \begin{equation}\label{eqn: cocycle}
        Z(\alpha\beta) = Z(\alpha)+\Ad_{\rho_0(\alpha)}Z(\beta),
    \end{equation}
    and
    \begin{equation}\label{eqn: 2nd order}
        W(\alpha\beta) = W(\alpha) + \Ad_{\rho_0(\alpha)}W(\beta) + Z(\alpha)\Ad_{\rho_0(\alpha)}Z(\beta). 
    \end{equation}
    Equations \eqref{eqn: traceless 1} and \eqref{eqn: cocycle} show that $Z \in Z^1(\pi,\mathfrak g_{\Ad_{\rho_0}})$, while \eqref{eqn: 2nd order} can be reformulated as
    \[\delta W(\alpha,\beta) = W(\alpha\beta) - W(\alpha) -\Ad_{\rho_0(\alpha)}W(\beta) = Z(\alpha)\Ad_{\rho_0(\alpha)}Z(\beta).\]
    
    We now compute the coboundary of $2W-Z^2$ (which is traceless by \eqref{eqn: traceless 2}).  The coboundary of $Z^2$ is
    \begin{align*}
    \delta Z^2(\alpha,\beta) &= (Z(\alpha\beta)^2) - Z(\alpha)^2 - (\Ad_{\rho_0(\alpha)}Z(\beta))^2 \\
    &= (Z(\alpha) + (\Ad_{\rho_0(\alpha)}Z(\beta)))^2 - Z(\alpha)^2 - (\Ad_{\rho_0(\alpha)}Z(\beta))^2 \\
    &= Z(\alpha)\Ad_{\rho_0(\alpha)}Z(\beta) + \Ad_{\rho_0(\alpha)}Z(\beta)Z(\alpha).
    \end{align*}
    Using linearity of $\delta$, we conclude that
    \begin{equation} \label{eqn: coboundary in H^2}
        \delta(2W - Z^2)(\alpha,\beta) = Z(\alpha)\Ad_{\rho_0(\alpha)}Z(\beta) - \Ad_{\rho_0(\alpha)}Z(\beta)Z(\alpha)  =[Z,Z](\alpha,\beta).  
    \end{equation}
    That is, $\delta(2W-Z^2) = [Z,Z]$.  
    
    Equations \eqref{eqn: coboundary in H^2} and \eqref{eqn: traceless 2} say that $[Z,Z]$ is a coboundary, which is to say that it vanishes in $H^2(\pi,\mathfrak g_{\Ad_{\rho_0}})$. 
\end{proof}

\begin{remark}\label{rmk: quadratic tangent}
    This identifies $Z^1(\pi,\mathfrak{g}_{\Ad_\rho})$ with the Zariski tangent space to the variety of homomorphisms from $\pi$ to $\SL(n,\K)$ at $\rho$ and the vanishing locus of the quadratic map 
    \[
    Z\in Z^1(\pi,\mathfrak{g}_{\Ad_\rho}) \mapsto [Z,Z]\in H^2(\pi, g_{\Ad_{\rho}})
    \]
    with its quadratic tangent cone at $\rho$; see \cite{Goldman:reps_of_surfaces}.
\end{remark}

\subsection{Deforming away from diagonal representations}
Let $\chi: \pi \to \K^\times$ be a non-trivial character, and let
\[\rho_\chi = \begin{pmatrix}
    \chi &  \\ & \chi\inverse
\end{pmatrix} \]
be the corresponding diagonal representation into $\SL(2,\K)$. The Lie algebra is denoted by $\mathfrak g$ (without a decoration indicating the field).  

Since 
\[\begin{pmatrix}
    \chi & \\& \chi\inverse
\end{pmatrix} \begin{pmatrix}
    x& y\\z&-x
\end{pmatrix}\begin{pmatrix}
    \chi\inverse&\\ & \chi
\end{pmatrix}=\begin{pmatrix}
    x& \chi^2y\\ \chi^{-2}z& -x
\end{pmatrix},\]
the coefficient module $\mathfrak g$ decomposes as a $\pi$-representation as 
\[\frak g_{\Ad_{\rho_\chi}}\cong \K_{\chi^2} \oplus \K_{\chi^{-2}}\oplus \K\]
inducing a decomposition in cohomology
\[H^1(\pi,\mathfrak g_{\Ad_{\rho_\chi}}) \cong H^1\left(\pi,\K_{\chi^{2}}\right) \oplus H^1\left(\pi,\K_{\chi^{-2}}\right) \oplus H^1(\pi,\K).\]
Since $\K_{
\chi^2}$ and $\K_{\chi^{-2}}$ are dual coefficient modules, the cup product and natural evaluation of coeffients defines an anti-symmetric bilinear pairing 
\[[\cdot \smile \cdot] : Z^1\left(\pi,\K_{\chi^{2}}\right) \times  Z^1\left(\pi,\K_{\chi^{-2}}\right) \to H^2(\pi, \K).\]
\begin{lemma}\label{lem: bracket and cup}
    Suppose $\rho: F_{2g} \to \SL(2,\K[t])$ factors through $\pi$ to second order, and suppose the reduction $\mod t$ is a non-trivial diagonal representation $\rho_\chi$.  Let $u \in Z^1(\pi, \K_{\chi^2})$ and $\ell \in Z^1(\pi, \K_{\chi^-2})$, and suppose that 
    \[\rho = (I + tZ +t^2W)\rho_\chi \mod t^3,\]
    where $Z = \begin{pmatrix}
        & u\\\ell&
    \end{pmatrix}$ and $W: \pi \to M_2(\mathbb K)$. Then $[u\smile \ell] = 0$.
\end{lemma}
\begin{proof}
    By computation, 
    \begin{align*}
        [Z,Z](\alpha,\beta) &= 
       X\left( u(\alpha)\chi^{-2}(\alpha)\ell(\beta) -\ell(\alpha)\chi^2(\alpha)u(\beta)\right) \\
       & = X(u\smile\ell - \ell\smile u)(\alpha,\beta),
    \end{align*}
    where $X = \begin{pmatrix}
        1& \\ & -1
    \end{pmatrix}$.
    In cohomology, $[u\smile \ell - \ell\smile u] = 2[u\smile \ell]$.
    By Proposition \ref{prop: deformation theory}, $[Z,Z]=0$.  We conclude then that $[u\smile \ell]=0$.
\end{proof}

\subsection{Zooming in at a diagonal representation}
We extract sub-limiting directions from a sequence $[\rho_n]$ converging to a non-trivial diagonal $[\rho_\chi]$.

In this section, $\mathcal X_{\K}$ is defined to be the quotient of the space of semi-simple representations into $\PSL(2,\K)$ modulo the $\PSL(2,\K)$-conjugation action. 
Let $\mathcal Z_\K\subset \mathcal X_\K$ be the locus of characters admitting a diagonal representative.  
It is known that the natural map $\mathcal X_{\R}\to \mathcal X$ is a homeomorphism (where $\mathcal X$ is defined as in the introduction); see, e.g., \cite[\S2.1]{Porti}.

The main idea in the following result is that we can find explicit $2$-jets of representations of the free group $F_{2g}$ starting at the diagonal locus and ending at $\rho_n$.  These $2$-jets subconverge by compactness, and any sublimit is a representation to second order.  An easy computation then furnishes estimates on traces.  

\begin{proposition}\label{prop: zooming in}
    Let $\chi \in H^1(\pi, \K^\times)$ be such that $\chi^2\not=1$ and let $[\rho_n] \notin \mathcal Z_\K$ be a sequence converging to $[\rho_\chi]$ in $\mathcal X_\K$.  
    Up to passing to a subsequence with the same name, there are $\epsilon_n\to 0$,  $u \in Z^1(\pi, \K_{\chi^2})$ and $\ell \in Z^1(\pi, \K_{\chi^{-2}})$, both different from zero in cohomology, such that 
    \begin{itemize}
        \item $[u\smile \ell] = 0 \in H^2(\pi,\K)$, and 
        \item for all $\eta \in [\pi,\pi]$, we have 
    \begin{equation} \label{eqn: hessian of trace}
        \frac{\tr(\rho_{n}(\eta))-2}{\epsilon_n^2}\to u(\eta)\ell(\eta), \text{ as } n \to \infty.
    \end{equation}
    \end{itemize}
\end{proposition}

\begin{remark}
    In fact, the quantity $u(\gamma)\ell(\gamma) \in \K$ only depends on the cohomology classes of $u$ and $\ell$.  See Proposition \ref{prop: sep curve bilinear form}.
\end{remark}
For the proof, it will be convenient to choose lifts to $\SL(2,\K)$, which is always possible in the Euler number zero component when $\K  = \R$ and in the component containing diagonal representations when $\K = \C$.
\begin{proof}
   Let $\gamma_1, ..., \gamma_{2g}$ be a standard generating set for $\pi$. 
   Without loss of generality, we may assume that $\chi(\gamma_1) \not=1$.  
   By definition of $\mathcal X_\K$ as the quotient of semi-simple representations modulo $\PSL(2,\K)$-conjugation, we may find lifts of $\rho_n$ to $\SL(2,\K)$ for which $\rho_n(\gamma_i) \to \rho(\gamma_i)$ for all $i$ as $n\to \infty$.  
   Thus, up to conjugation in $\SL(2,\K)$ by elements tending to $I$, we may assume that for $n$ large enough 
   \[\rho_n(\gamma_1) = \begin{pmatrix}
       \chi_n & \\ & \chi_n\inverse
   \end{pmatrix},\]
    where $\chi_n\to \chi(\gamma_1)$ as $n\to \infty$.  
   We write 
   \[\rho_n(\gamma_i) = \begin{pmatrix}
       a_{i,n} & b_{i,n}\\c_{i,n} & d_{i,n}
   \end{pmatrix},\]
   where all entries are in the field $\K$.
    Let $b_n = (b_{1,n}, \ldots b_{2g,n}) \in \K^{2g}$, and $c_n = (c_{1,n}, \ldots, c_{2g,n}) \in \K^{2g}$.  
    Conjugating $\rho_n$ by $\diag(t,t\inverse)$ scales $b_n$ by $t^2$ and $c_n$ by $t^{-2}$.  Using the usual norm $\|\cdot \|$  on $\K^{2g}$, we may thus assume that 
    \[\|b_n\| = \|c_n\| = \epsilon_n>0.\]
    Then $\epsilon_n \to 0$.  Indeed, the products $b_{i,n}c_{j,n}$ are invariant by $\diag(t,t\inverse)$ conjugation for all $i$ and $j$, and we know that $b_{i,n}$ and $c_{i,n}$ both tend to $0$ as $n \to \infty$ for all $i$.  Hence, $\epsilon_n^4 = \|b_n\|^2\|c_n\|^2 = \sum_{i,j}|b_{i,n}|^2|c_{j,n}|^2 \to 0$, as $n \to \infty$. 

    Write $b_{i,n} = \epsilon_n\beta_{i,n}$ and $c_{i,n} = \epsilon_n\zeta_{i,n}$.  
    Thus, up to taking a subsequence (without renaming), we can assume that 
    \[b_n/\epsilon_n \to \beta= (\beta_1,\ldots, \beta_{2g}) \in \K^{2g} \text{ and } c_n/\epsilon_n\to \zeta = (\zeta_1, \ldots, \zeta_{2g}) \in \K^{2g}.\]
    Then both $\beta$ and $\zeta$ are non-zero, because they have unit norm.
    \medskip

    For each $n$ and $i$, define
    \[M_{i,n} (t) = \begin{pmatrix}
        a_{i,n} & t\beta_{i,n} \\ t\zeta_{i,n} & \dfrac{1+t^2\beta_{i,n}\zeta_{i,n}}{a_{i,n}}
    \end{pmatrix}.\]
    The path $\rho_n^t$ defined by $\rho_n^t(\gamma_i) = M_{i,n}(t)$ is a representation of the free group $F_{2g}$ generated by the $\gamma_i$'s into $\SL(2,\K)$.  However, $\rho_n^0$ is diagonal, hence has abelian image and factors through $\pi$.  Moreover, since $a_{i,n}d_{i,n} - \epsilon_n^2 \beta_{i,n}\zeta_{i,n} = 1$, we see that $\rho_n^{\epsilon_n} = \rho_n$, and so also factors through $\pi$.

    As $n\to \infty$, we have 
    \[M_{i,n}(t) \to M_i(t) = \begin{pmatrix}
        \chi_i & t\beta_i\\ 
        t\zeta_i & \dfrac{1+t^2\beta_i\zeta_i}{\chi_i}
    \end{pmatrix},\]
    where $\chi_i = \chi(\gamma_i)$.
    We again define a path $\rho^t$ of representations of the free group $F_{2g}$ into $\SL(2,\K[t])$ by $\rho^t(\gamma_i) = M_i(t)$.  Observe that $\rho^0 =\rho_\chi$.  

    \begin{claim}\label{clm: relation}
        The path $\rho^t$ is a representation of $\pi$, up to second order, i.e.,
        \[\prod_{i = 1}^g[\rho^t(\gamma_{2i-1}),\rho^t(\gamma_{2i})]= I \mod t^3.\]
    \end{claim}

    \begin{proof}[Proof of the claim]
        Let $r = \prod_{i = 1}^g [\gamma_{2i-1},\gamma_{2i}]$, and define $R_n(t) = \rho_n^t(r)$.  Then $R_n(0)= I$, so expanding about $t=0$, we get
        \[R_n(t) = I +tA_n +t^2B_n+O(t^3),\]
        where the implicit constants are uniform in $n$, because the matrices  $\rho_n^t(\gamma_i)$ converge to $\rho^t(\gamma_i)$ for all $i$. 
        Similarly, for $R(t) = \rho^t(r)$, we have 
        \[R(t) = I + tA+t^2B +O(t^3),\]
        where $A_n \to A$ and $B_n\to B$ in $M_2(\K)$ as $n \to \infty$.

        The diagonal entries of $M_{i,n}(t)$ are even functions of $t$, while the off diagonal entries are odd functions of $t$.  This property is preserved by inversion and multiplication in $\SL(2,\K[t])$.  Therefore, $A_n$ is an off-diagonal matrix and $B_n$ is diagonal.  

        Since $R_n(\epsilon_n) = I$, we get 
        \[0 = \epsilon_nA_n + O(\epsilon_n^3) \text{ and } 0 = \epsilon_n^2B_n +O(\epsilon_n^4).\]
        In particular, $A_n = O(\epsilon_n^2)$ and $B_n = O(\epsilon_n^2)$.  Since $\epsilon_n \to 0$, we obtain $A = B = 0$, i.e., $R(t) = I \mod t^3$.
    \end{proof}

    The claim allows us to consider $\rho^t$ as a representation $\pi \to \SL(2,\K[t]/(t^3))$.  As such 
    \[\rho^t(\gamma) = \left( I +tZ(\gamma)+t^2W(\gamma)\right)\rho_\chi(\gamma),\]
    for all $\gamma \in \pi$.
    Since $Z(\gamma)$ is only non-zero in its off-diagonal entries, Proposition \ref{prop: deformation theory} and Lemma \ref{lem: bracket and cup} give that 
    \[Z(\gamma) = \begin{pmatrix}
        & u(\gamma)\\\ell(\gamma)&
    \end{pmatrix},\]
    where $u \in Z^1(\pi, \K_{\chi^2})$ and $\ell \in Z^1(\pi, \K_{\chi^{-2}})$ satisfy $[u\smile \ell]=0$.  Furthermore, from the definitions, it follows that $u(\gamma_i) = \chi_i\beta_i$ and $\ell(\gamma_i) = \chi_i\inverse \zeta_i$.
    \medskip

    To see that $u$ and $\ell$ are not coboundaries, observe that $u(\gamma_1) = \chi_1\beta_1 =0$ and $\ell(\gamma_1) = \chi_1\inverse \zeta_1 = 0$, because we chose $\rho_n$ so that $\rho_n(\gamma_1)$ was diagonal.  Since $\chi(\gamma_1)\not=1$ (by assumption), it follows that a coboundary that vanishes on $\gamma_1$ must vanish everywhere, since the coboundaries are of the form $\gamma \mapsto (1-\chi^2(\gamma))x$ or $\gamma \mapsto (1-\chi^{-2}(\gamma))x$ for $x \in \K$.  Since $u$ and $\ell$ are non-zero, they are thus not coboundaries.  
    \medskip

    Finally, we need to verify the assertion about the behavior of traces of elements of $[\pi,\pi]$.
    For any $\gamma \in \pi$ 
    the diagonal entries of $\rho_n^t(\gamma)$ are even functions of $t$ and the off-diagonal entries are odd functions of $t$.   For $\eta \in [\pi,\pi]$, we have $\rho_n^0$ is diagonal so $\rho_n^0(\eta) = I$, and
    \[\rho_n^t(\eta) = \begin{pmatrix}
        1 +p_n(t) & q_n(t)\\r_n(t) & 1+s_n(t)
    \end{pmatrix},\]
    where $p_n(t)$, $ s_n(t) = O(t^2)$, while $q_n(t)$, $r_n(t) = O(t)$.  Convergence of $\rho_n^t$ to $\rho^t$ gives 
    \[q_n(\epsilon_n)/\epsilon_n\to u(\eta)\chi(\eta)\inverse = u(\eta)\]
    because $\eta \in [\pi,\pi]$.  Similarly, \[r_n(\epsilon_n)/\epsilon_n \to \ell(\eta) \chi(\eta) = \ell(\eta).\]

    Now, using $\det (\rho_n^t(\eta))=1$, we see that \[p_n +s_n =  q_nr_n -p_ns_n.\]  
    At $t = \epsilon_n$, the left hand side is $\tr(\rho_n(\eta))-2$.
    Since $p_n(\epsilon_n) s_n(\epsilon_n) = O(\epsilon_n^4)$, we conclude that
    \[\frac{\tr(\rho_n(\eta))-2}{\epsilon_n^2} = \frac{q_n(\epsilon_n)r_n(\epsilon_n)}{\epsilon_n^2}+O(\epsilon_n^2) \to u(\eta)\ell(\eta),  \]
    as claimed.
\end{proof}

\section{Separating curves and intersection computations}

Recall that a simple curve which lies in the commutator subgroup $[\pi,\pi]$ is a separating curve. Proposition \ref{prop: zooming in} gives a criterion for showing that a given separating curve $\eta$ is eventually elliptic for a given sequence of representations converging to a diagonal representation; namely, if $u(\eta)\ell(\eta) < 0$ for the limiting pair of cocycles $(u,\ell)$.  To prove our Main Theorem we need to show that for any sequence of representations converging to a diagonal representation, some separating curve $\eta$ satisfies this criterion. We will do so by showing that for each separating curve $\eta$, there is a nonempty open set of the space of pairs $(u,\ell)$ for which $u(\eta)\ell(\eta) < 0$ and then showing that the orbit of this open set under the action of the mapping class group covers the space of all possible pairs. 

In this section we give an alternative description of the space of pairs $(u,\ell)$ in terms of local systems, which will allow us to make explicit computation of the quantity $u(\eta)\ell(\eta)$ by topological arguments. The first subsection introduces background on cohomology and homology with coefficients in a local system. The second subsection contains the main topological argument.

\subsection{Chains, cochains, and intersection for line bundles}\label{subsec: chains and intersection}

Consider the flat line bundle $E_\chi$ over $S$ with holonomy representation $\chi: \pi\to \R_{>0}$, i.e., the fibers of $E_\chi$ are copies of $\R$, viewed as a one-dimensional vector space.\footnote{When comparing this section to the rest of the paper, one should make the substitution $\chi = \chi^2$, throughout.}
We discuss briefly the homology and cohomology groups $H_1(E_\chi)$ and $H^1(E_\chi)$.
An introduction to these objects can be found in \cite[\S4]{JM:deformations}.

Let $[0,\ldots, i]$ denote the standard, oriented $i$-simplex with vertices labeled $0, ..., i$.
An $i$-chain is a flat section over a smooth singular $i$-simplex $\sigma: [0,...,i]\to S$, which is determined by its value in the fiber over $\sigma(0)$.
The boundary maps from the vector space of $i$-chains to the $(i-1)$-chains are defined, as usual, taking into account the flat connection on $E_\chi$. 
The cycles are denoted $Z_i(E_\chi)$, the boundaries are $B_i(E_\chi)$, and the quotient homology groups are $H_i(E_\chi)$.

Any closed curve $\gamma\subset S$ with $\chi(\gamma) =1$, e.g., any separating curve, and any non-zero vector $v$ in the fiber in $E_\chi$ over a point of $\gamma$ determines 
a flat section of $E_\chi\vert_\gamma$,
which defines an element $[\gamma_v]\in H_1(E_\chi)$.

Consider also the dual flat vector bundle $E_{\chi\inverse}$, which we identify with $E_\chi^*$.
There is a natural anti-symmetric bilinear \emph{intersection pairing} 
\[i: H_1(E_\chi)\times H_1(E_{\chi\inverse})\to \R.\]
See \cite[Lemma 4.3]{JM:deformations}. For transverse closed curves $\alpha$ and $\beta$ with $\chi(\alpha) = \chi(\beta) = 1$ and flat sections $\alpha_v$ and $\beta_w$ of $E_\chi$ and $E_{\chi\inverse}$, respectively,
we have 
\[i([\alpha_v], [\beta_w]) = \sum_{q \in \alpha\pitchfork \beta} \epsilon(q, \alpha, \beta) \alpha_v(q)\beta_w(q)  \in \R,\]
where $\epsilon(q,\alpha,\beta)$ is the sign of the intersection between $\alpha$ and $\beta$ at $q$.  

In particular, for any simple separating curve $\eta\subset S$, and $[\eta_v] \in H_1(E_\chi)$ and $[\eta_w] \in H_1(E_{\chi\inverse})$, we have 
\[i([\eta_v], [\eta_w])=0.\]
This statement contains two claims which are easily verified: first that the necessary sections over $\gamma$ always exist, and second that the intersection of a simple curve with itself is trivial.

\medskip

Now we consider the cohomology group $H^1(E_\chi)$.
An $i$-cochain $z$ assigns to every smooth singular $i$-simplex $\sigma: [0,\ldots, i]\to S$ an element of the fiber of $E_\chi$ over $\sigma(0)$.
The coboundary operators are again defined as usual, the $i$-cocycles are denoted $Z^i(E_\chi)$ and the cohomology groups are $H^i(E_\chi)$.
The cocycle condition when $i=1$ is as follows.
Whenever $\gamma_1$, $\gamma_2$, and $\gamma_3$ are $1$-simplices on $S$ with $\gamma_1\cdot \gamma_2 \sim \gamma_3$ rel endpoints, 
we have
\[z(\gamma_3 ) = z(\gamma_1) + \overline {\gamma_1}_* z(\gamma_2),\]
where $\overline {\gamma_1}_*$ is the parallel transport taking the value of $z(\gamma_2)$ over $\gamma_2(0)$ to $\gamma_1(0)$, and the equality is interpreted in the fiber over $\gamma_1(0) = \gamma_3(0)$. 

There is a natural \emph{Poincar\'e duality} linear isomorphism 
\[PD: H_1(E_\chi) \to H^1(E_\chi),\]
described at the level of cochains as follows.
For the class $[\gamma_v] \in H_1(E_\chi)$ described as above, 
$PD(\gamma_v)(\sigma)$ is defined for transverse cycles by
\[PD(\gamma_v)(\sigma) = 
\sum_{q \in \gamma \pitchfork \sigma} \epsilon(q, \gamma, \sigma) \overline{\sigma_*}\gamma_v(q) \in E_\chi\vert_{\sigma(0)}\]
In other words, this is the flat section over $\sigma(0)$ defined by summing together the value of $\gamma_v$ at all points of intersection, with a sign from the orientation of intersection at each point. Note that we may also think of a vector in the fiber over $\sigma(0)$ as defining a flat section over $\sigma$.

For $[\gamma_v] \in H_1(E_\chi)$ as above, we use the notation
\[\iota_{\gamma_v}^\chi :=PD([\gamma_v])\in H^1(E_\chi),\]
or just $\iota_\gamma^\chi$, when the choice of $v$ is immaterial.  
We remark that  $\iota_{\alpha_v}^\chi \smile\iota^{\chi\inverse}_{\beta_w} \in H^2(E_\chi\otimes E_{\chi\inverse})\cong H^2(S,\R)\cong \R$, and
\[\int_S\iota_{\alpha_v}^\chi \smile\iota^{\chi\inverse}_{\beta_w} = i([\alpha_v], [\beta_w]), \]
where we use the integral over $S$ to mean the evaluation on the fundamental class of $S$.

Evaluating a cocycle on based loops in $\pi = \pi_1(S,x_0)$ defines a natural linear map
\[H^1(E_\chi) \to H^1(\pi,\R_\chi),\]
which is an isomorphism \cite{JM:deformations}.

\subsection{Nonvanishing of the intersection pairing on separating curves}\label{subsec: computations}

By a separating curve $\eta \subset S$, we mean a closed curve whose image in $H_1(S,\R)$ is trivial.  Homotopically essential simple closed curves that separate $S$ are examples.  

\begin{proposition}\label{prop: sep curve bilinear form}
Every separating curve $\eta$ gives a well-defined bilinear form
    \[
    B_\eta: H^1(E_\chi) \times H^1(E_{\chi^{-1}}) \to \R
    \]
by $B_\eta([u],[\ell]) = \langle u(\eta),\ell(\eta)\rangle$, where the pairing takes place in any fiber over $\eta$. $B_\eta$ only depends on the homotopy class of $\eta$, and is independent of its orientation.
\end{proposition}

\begin{proof}
    The fiberwise duality between $E_\chi$ and $E_{\chi^{-1}}$ makes it clear that the pairing between the sections $u(\eta)$ and $\ell(\eta)$ is constant on $\eta$, and is invariant under homotopy of $\eta$. It just remains to check that $B_\eta$ vanishes when $u$ or $\ell$ is a coboundary.  Choose a parameterization of $\eta$ by the interval $[0,1]$.  Suppose without loss of generality that $u = \delta z$, where $z$ is a section of $E_\chi$. Then $u(\eta) = \overline{\eta}_* z(\eta(1)) - z(\eta(0))$. This is equal to zero since $\eta$ is separating and so the holonomy of $E_\chi$ over $\eta$ is trivial.  

    Since $u$ is a cocycle, $u(\overline{\eta}) = - \eta_* u(\eta)$. Since the holonomy over eta is trivial, this is equal to $-u(\eta)$, and similarly for $\ell$. Therefore reversing the orientation of $\eta$ introduces two minus signs, which cancel out.
\end{proof}

\begin{remark}
    Under the identification between $H^1(E_\chi)$ and $H^1(\pi, \R_\chi)$, the bilinear form $B_\eta$ becomes the pairing $u(\eta)\ell(\eta)$ of \eqref{eqn: hessian of trace} in Proposition \ref{prop: zooming in}. Hence Proposition \ref{prop: zooming in} implies that the natural inclusion 
    \[
    Z^1(\pi, \R_{\chi^2}\oplus \R_{\chi^{-2}}) \to Z^1(\pi, \mathfrak{g}_{\Ad_{\rho_\chi}})
    \]
    identifies $B_\eta$ with the restriction of the Hessian of $\tr_\eta$ on the variety of homomorphisms (see Remark \ref{rmk: quadratic tangent}).
\end{remark}

 Let $V = H^1(E_\chi)$, and recall that Poincar\'e duality and the intersection pairing identifies $V^*$ with $H^1(E_{\chi\inverse})$.
 Then $H^1(E_\chi)\otimes H^1(E_{\chi\inverse}) \cong \End(V)$, and this identification gives 
 \[\tr([u]\otimes [\ell]) = [u\smile \ell]\in H^2(S,\R)\cong \R.\]

\begin{figure}[h]
    \centering
    \includegraphics[width=0.4\linewidth]{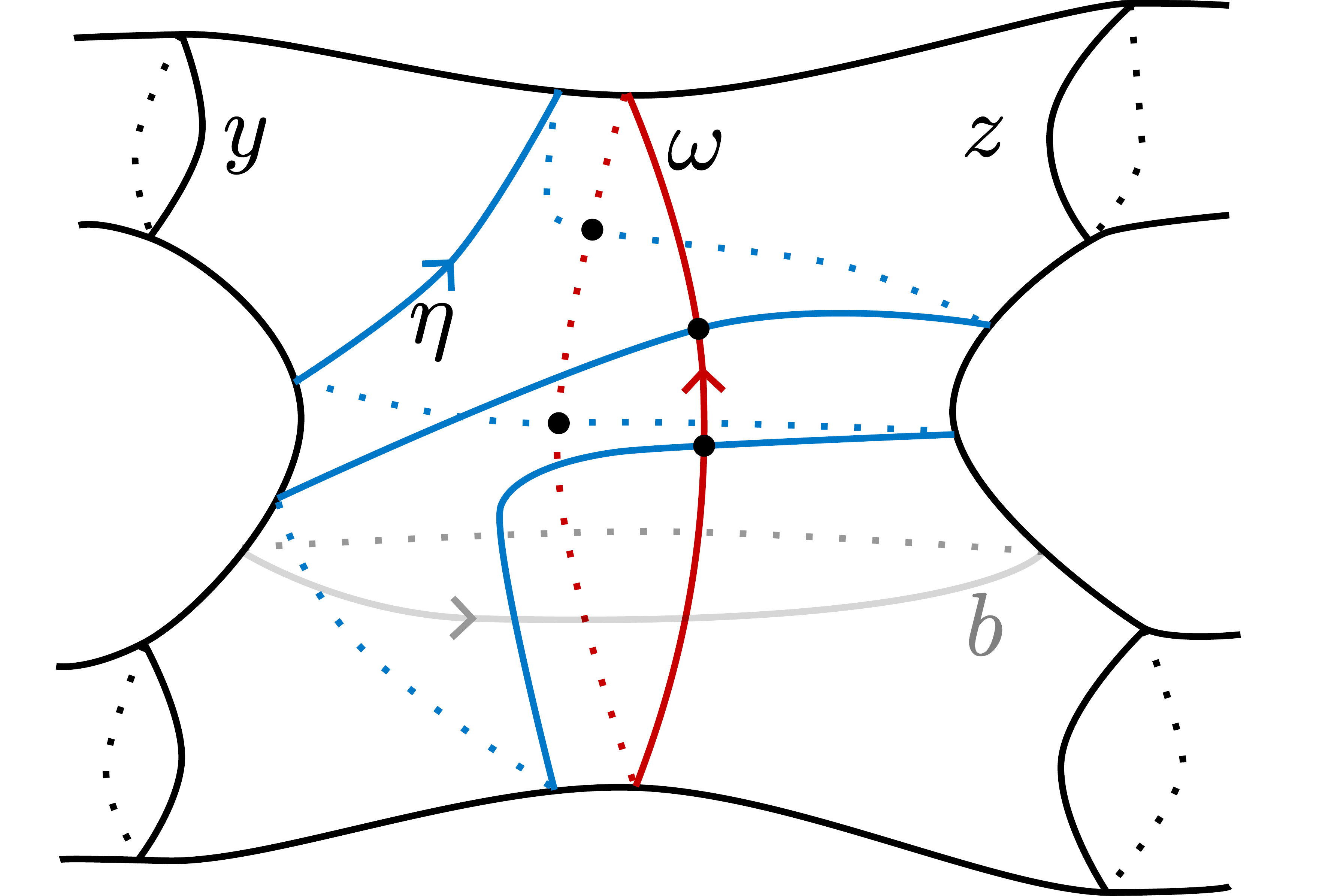}
    \caption{The subsurface $S'\subset S$ and the curves $\omega$ and $\eta$ involved in the computation in this section. They differ by a Dehn twist about the curve $b$. Their four points of intersection are shown.}
    \label{fig: comp1}
\end{figure}

\begin{figure}[h]
    \centering
    \includegraphics[width=0.4\linewidth]{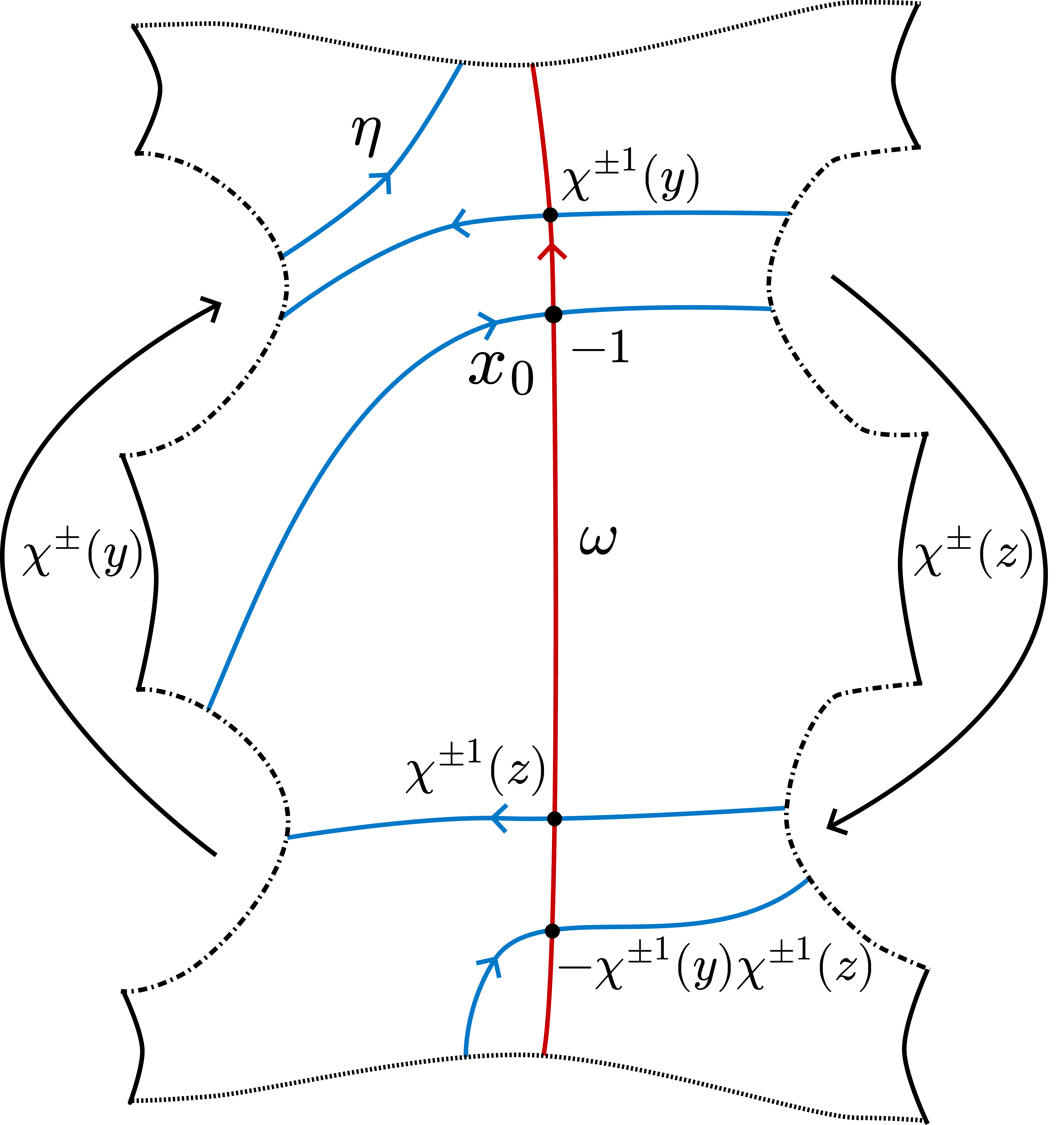}
    \caption{The same two curves $\omega$ and $\eta$ from Figure \ref{fig: comp1} with the surface unfolded. The curves $y$ and $z$ are gluing where indicated. To compute intersection in $E_{\chi^{\pm 1}}$, we trivialize the bundle over the disk shown, and push the holonomy onto the gluing loops. The four intersection points each contribute the indicated term to the intersection $\iota^{\chi^{\pm 1}}_\omega(\eta)$.}
    \label{fig: comp2}
\end{figure}

Now, we produce a very concrete rank one simple tensor to pin down how it evaluates on a specific separating curve. 
\begin{lemma}\label{lem: some chi have a good direction}
    With curves as in Figure \ref{fig: comp1}, for all $\chi:\pi\to \R_{>0}$ satisfying
    \[(\chi(y)-1)(\chi(z)-1)\neq 0 \]
    there exists a rank one traceless element 
    \[[u] \otimes [\ell] \in H^1(E_\chi)\otimes H^1(E_{\chi\inverse})\]
    such that
    \[ B_\eta([u],[\ell]) <0.\]

\end{lemma}

\begin{proof}
    We choose a pair of separating curves $\omega$ and $\eta$, arranged in a subsurface $S'\subset S$ as in Figure \ref{fig: comp1}.
In particular, there exists a curve $b$ which intersects $\omega$ twice (geometrically) and $\eta$ differs from $\omega$ by a Dehn twist around $b$. Such a pair always exists when genus is greater than one.

To derive equations \eqref{eqn: see fig 2} and \eqref{eqn: see fig 2 again} we compute using Figure \ref{fig: comp2} following the definitions in \S\ref{subsec: chains and intersection}. 
Choose $\iota_\omega^\chi$ and $\iota_\omega^{\chi\inverse}$ to be Poincar\'e dual to a pair of sections of $E_\chi$ and $E_{\chi\inverse}$ over $\omega$ that pair to $1$ pointwise. 
Choosing the basepoint $x_0$ to trivialize the bundle $E_\chi$ (or $E_{\chi\inverse}$), we find
\begin{equation} \label{eqn: see fig 2}
    \iota^\chi_\omega(\eta) = -1+\chi(y)+\chi(z)-\chi(y)\chi(z)=-(\chi(y)-1)(\chi(z)-1)
\end{equation}
and similarly
\begin{equation} \label{eqn: see fig 2 again}
   \iota^{\chi\inverse}_\omega(\eta) = -(\chi\inverse(y)-1)(\chi\inverse(z)-1). 
\end{equation}

    Let $[u] \otimes [\ell] = - \iota^\chi_\omega\otimes\iota^{\chi\inverse}_\omega$, so that 
    \[ B_\eta([u],[\ell]) = - \frac{(\chi(y)-1)^2(\chi(z)-1)^2}{\chi(y)\chi(z)}<0.\]
    The traceless condition is exactly that $\iota_\omega^\chi \smile \iota_\omega^{\chi\inverse} = i([\omega], [\omega])= 0$.
    This completes the proof.
\end{proof}

\begin{remark}\label{remark: scale the intersection}
    In the computation we have suppressed the vector defining the section over $\omega$ (it is defined by an implicit identification $ E_{\chi^{\pm1}}\vert_{x_0}\cong \R$). Changing this choice has the effect of scaling the intersections evaluated above.
\end{remark}


We obtain the following via the action of the mapping class group on homology.

\begin{proposition}\label{prop: all chi have a good direction}
    For all non-trivial $\chi:\pi\to \R_{>0}$, there exists a rank one, traceless element 
    \[[u]\otimes[\ell]\in H^1(E_\chi)\otimes H^1(E_{\chi\inverse})\]
    and simple separating curve $\eta$ such that 
    \[B_\eta([u],[\ell]) <0.\]
\end{proposition}

\begin{proof}
    Observe that $\chi=e^\varphi$ for some $\varphi\in H^1(S,\R)$.
    With $y$ and $z$ defined as in Figure \ref{fig: comp2}, define $U = \{\varphi\in H^1(S,\R) : \varphi(y) \varphi(z) \not=0\}$.  This set is a cone, i.e., invariant under scaling by $\R$.  
    For every $\varphi \in U$, the character $e^{\varphi}$ satisfies the hypotheses of Lemma \ref{lem: some chi have a good direction}.
    
    It is well-known that the action of the mapping class group on $H^1(S,\R)$ (which is identified with the action of $\mathrm{Sp}(2g,\Z)$) acts minimally on the projective space $\mathbb P(H^1(S,\R))$.
    Since $U$ is a cone, minimality implies that the orbit of every nonzero $\varphi$ meets $U$. Applying Lemma \ref{lem: some chi have a good direction} there and transporting the tensor and curve back proves the proposition.
\end{proof}

The following result of Ghazouani is part of the main theorem of \cite{Selim:orbitclosure}. 

\begin{theorem}[{\cite[Theorem 5.1 (1)]{Selim:orbitclosure}}]\label{thm: G}
    Let $\chi: \pi \to \R_{>0}$ be a non-trivial character.   Then the stabilizer $\textrm{stab}(\chi)<\Mod S$ acts naturally on the projectivization of
    $H^1(\pi,\R_\chi)$ by projective transformations and has image which is dense in $\PSL(H^1(S,\R_\chi))\cong \PSL(\R^{2g-2})$.
\end{theorem}
\begin{remark}
    Ghazouani's theorem states that the image is dense in $\PGL(\R^{2g-2})$, but the proof gives density in $\PSL(\R^{2g-2})$.  This distinction is important; see Remark \ref{rmk: genus 2}. 
\end{remark}

\begin{remark}
    Recall that $\Mod(S)$ maps naturally to $\Out(\pi)$.
    The stabilizer of $\chi$ in $\Aut(\pi)$ acts linearly  on $H^1(\pi,\R_\chi)$ and $H^1(\pi,\R_{\chi\inverse})$, but conjugation by an element $\gamma \in \pi$ effects multiplication by $\chi(\gamma)$ on $H^1(\pi, \R_\chi)$ and $\chi\inverse(\gamma)$ on $H^1(\pi, \R_{\chi\inverse})$.  Taking the tensor product, the linear action of the stabilizer of $\chi$ in $\Aut(\pi)$ thus descends to a linear action of the stabilizer of $\chi$ in $\Out(\pi)$.  Theorem \ref{thm: G} thus gives that $\textrm{stab}(\chi)<\Mod(S)$ acts by linear transformations on $H^1(\pi, \R_\chi)\otimes H^1(\pi,\R_{\chi\inverse})$, and that the closure of the image of $\textrm{stab}(\chi)$ in $\GL(H^1(\pi, \R_\chi)\otimes H^1(\pi, \R_{\chi\inverse}))$ is a copy of $\PSL(2g-2,\R)$.
\end{remark}

To prove this, one checks that different automorphisms representing a given mapping class differ by central (scalar) elements of $\GL(H^1(\pi,\R_\chi))$. The center acts trivially on the tensor product, as the two spaces are dual (for any vector space, the action of $Z(\GL(V))$ on $V\otimes V^*$ is trivial).

Combining Ghazouani's result with the previous corollary, we obtain 

\begin{theorem}\label{thm: negative pairing}
    Suppose the genus of $S$ is at least $3$.  Then for all non-trivial characters $\chi : \pi \to \R_{>0}$ and all $[u] \otimes [\ell] \in H^1(E_\chi)\otimes H^1(E_{\chi\inverse})$ which are rank one and traceless, there exists a simple separating curve $\eta$ such that 
    \[B_\eta([u],[\ell]) <0.\]
\end{theorem}

\begin{proof}
    Proposition \ref{prop: all chi have a good direction} provides an element $[u]\otimes[\ell]\in H^1(E_\chi)\otimes H^1(E_{\chi\inverse})$ and separating curve $\eta$ for which 
    the statement is satisfied. By continuity of the evaluation 
    \[B_\eta([u],[\ell]),\]
    there exists an open positive cone in $H^1(E_\chi)\otimes H^1(E_{\chi\inverse})$ such that the evaluation with any rank one element in this cone remains negative.
    When the genus is at least 3, $\PSL(2g-2,\R)$ acts transitively on the cone of traceless rank-one elements in $\End(\R^{2g-2})$. Hence, Ghazouani's density result (Theorem \ref{thm: G}) implies that every traceless rank-one element is in the mapping class group orbit of this cone, and hence pairs negatively with some separating curve.
\end{proof}

\begin{remark}\label{rmk: genus 2}
    Theorem \ref{thm: negative pairing} fails for genus $2$, because when $\dim H^1(\pi, \R_\chi) = 2g-2 = 2$, $\PSL(2,\R)$ has two distinct orbits on the set of traceless rank one elements of $\R^2\otimes (\R^2)^*$. Indeed, this action is the same as the action of the connected component of $SO(2,1)$ on the light cone in Minkowski 3-space $\R^{2,1}$). $\PSL(2,\R)$ preserves each component of the light cone, while $\PGL(2,\R)$ swaps them. 
    
    In fact, by \cite{MW:genustwo}, there are two open $\Mod(S)$-invariant sets in $\mathcal X_0$ for genus $2$.  
\end{remark}

\section{Proof of Theorem \ref{thm: exists neighborhood R}}

Suppose $[\rho_n] \in \mathcal X \setminus \mathcal Z$ is a sequence converging to a non-trivial diagonal representation $[\rho_\chi]$.  Passing to a subsequence without renaming, Proposition \ref{prop: zooming in} gives us numbers $\epsilon_n \to 0$, and cocycles $u \in Z^1(\pi, \R_{\chi^2})$ and $\ell \in Z^1(\pi, \R)_{\chi^{-2}}$ with $[u]\smile [\ell] = 0$ such that for all curves $\eta \in \pi$ whose corresponding free homotopy class is a simple closed separating curve, 
\[\lim_{n\to \infty}\frac{\tr(\rho_n(\eta))-2}{\epsilon_n^2} = u(\eta)\ell(\eta).\]
Then $[u]\otimes [\ell]$ is rank-one and traceless, because $[u\smile \ell] =0$.
Applying Theorem \ref{thm: negative pairing} to $u\otimes \ell$ supplies a  separating curve $\eta$ with $u(\eta)\ell(\eta)<0$.  
Thus, 
\[\tr(\rho_n(\eta))<2\]
for $n$ large enough.
Since $\tr(\rho_n(\eta)) \to 2$ as $n\to \infty$ (because $\rho_\chi(\eta) = I$), we conclude that $\rho_n(\eta)$ is elliptic for $n$ large enough.

We have proved that the set
\[V= \{[\rho] \in \mathcal X: \exists \text{ simple separating curve } \eta ~s.t. ~|\tr(\rho(\eta))|\le 2\}\]
contains an open neighborhood $U$ of $\mathcal Z\setminus \{[id]\}$, because $\mathcal X$ is first countable and every sequence of points converging to a point of $\mathcal Z\setminus \{[id]\}$ is eventually contained in $V$.  Furthermore, we proved that if $[\rho] \in V\setminus \mathcal Z$ is close enough to $\mathcal Z \setminus \{[id]\}$, then also there is a simple separating curve $\eta$ for which $|\tr(\rho(\eta))|<2$.
Thus the theorem is proved.

\bibliography{biblio}{}
\bibliographystyle{amsalpha.bst}
\end{document}